\documentclass{article}
\usepackage{amsthm,amsmath,amssymb,amscd, mathrsfs}
\usepackage[latin1]{inputenc}
\usepackage[all]{xy}

\numberwithin{equation}{section}

 \newtheorem{thm}[equation]{Theorem}
 \newtheorem{prop}[equation]{Proposition}

 \theoremstyle{definition}
 \newtheorem{definition}[equation]{Definition}
 \newtheorem{example}[equation]{Example}

 \theoremstyle{remark}
 \newtheorem{remark}[equation]{Remark}

 \def\Spec{{\rm Spec~}}
 \def\Spf{{\rm Spf}}
 \def\FMod{{\rm FMod}}

\def\Gr{{\rm Gr}}
 \def\Nilp{{\rm Nilp}}
 \def\Sets{{\rm Sets}}

 \def\an{{\rm an}}
 \def\et{\text{\rm \'et}}

\def\an{{\rm an}}
 
 \def\dR{{\rm dR}}

 \def\det{{\rm det}}

 \def\dom{{\rm dom}}
\def\an{{\rm an}}

\def\GL{{\rm GL}}
\def\Rep{{\rm Rep}}

\def\alg{{\rm alg}}

\def\CG{\mathcal{G}}
\def\Res{{\rm Res}}
\def\Hom{{\rm Hom}}
\def\Nilp{{\rm Nilp}}
\def\BF{\mathbb{F}}
\def\Flag{{\rm Flag}}
\def\GG{{\mathfrak G}}
\def\G{\mathbf G}
\def\Geta{G}
\def\dbl{[\![}
\def\dbr{]\!]}
\def\dpl{(\!(}
\def\dpr{)\!)}

\def\isoto{\stackrel{}{\mbox{\hspace{1mm}\raisebox{+1.4mm}{$\sim$}\hspace{-3.5mm}$\longrightarrow$}}}

\def\NilpFq{{\Nilp_{\mathbb{F}_{q}\dbl \zeta\dbr}}}

\begin{document}
\begin{title}
{Moduli spaces of local $\G$-shtukas}
\end{title}
\author{Eva Viehmann}
\date{}

\maketitle

Math subject classification 2010:
11G09,  
(11G18,  
14L05,  
14M15) 

\begin{abstract}{We give an overview of the theory of local $\G$-shtukas and their moduli spaces that were introduced in joint work of U.~Hartl and the author, and in the past years studied by many people. We also discuss relations to moduli of global $\G$-shtukas, properties of their special fiber through affine Deligne-Lusztig varieties and of their generic fiber, such as the period map.}
\end{abstract}
\section{Introduction}\label{secintro}
Local $\G$-shtukas are an analog over local function fields of $p$-divisible groups with additional structure. In this article we give an overview about the theory of moduli spaces of local $\G$-shtukas and their relation to moduli of global $\GG$-shtukas. It parallels the theory of Rapoport-Zink moduli spaces of $p$-divisible groups and their relation to Shimura varieties. Yet it has the additional charm that additional structure given by any parahoric group scheme $\G$ can be easily encoded and treated in a group-theoretic way.

We begin by defining local $\G$-shtukas for a parahoric group scheme $\G$ over $\Spec\BF_q\dbl z\dbr$ (see Section \ref{sec2}). They are pairs consisting of an $L^+{\G}$-torsor $\mathcal G$ together with an isomorphism $\sigma^*L\mathcal{G}\rightarrow L\mathcal{G}$ for the associated $L\G$-torsor. Here for unexplained notation we refer to the respective sections. Local $\G$-shtukas are the function field analog of $p$-divisible groups with additional structure, as well as a group theoretic generalization of Drinfeld's shtukas. 
Next we discuss possible bounds on the singularities of local $\G$-shtukas. These are a more general replacement of the minuscule coweights fixed to define Shimura varieties. In Section \ref{sec3} we consider deformations for local $\G$-shtukas, and the description of the universal deformation in terms of the loop group $L\G$. In Section \ref{sec4} we describe the analog of Rapoport-Zink moduli spaces for local $\G$-shtukas bounded by a given bound $\hat Z$ and within a given quasi-isogeny class. They are representable by a formal scheme, locally formally of finite type over $\Spf \breve R$. Here, $\breve R$ is the completion of the maximal unramified extension of the ring of definition of the chosen bound $\hat Z$. Also, they admit a tower of coverings indexed by compact open subgroups of $\G(\BF_q\dpl z\dpr)$. Our next topic is the special fiber of a Rapoport-Zink moduli space of local $\G$-shtukas. It can be identified with a so-called affine Deligne-Lusztig variety. In Section \ref{sec6} we discuss applications of this result to geometric questions such as dimensions and closure relations of Newton strata. In Section \ref{sec5} we describe the relation between our local $\G$-shtukas and global $\GG$-shtukas where $\GG$ is a parahoric group scheme over a smooth projective geometrically irreducible curve $C$ over $\mathbb{F}_q$. There are several results in this direction, the central idea being that one wants to describe a global $\GG$-shtuka by corresponding local $\G_i$-shtukas associated with each of the legs $c_i$. Due to our choice of definition of local $\G$-shtukas, this works particularly well when considering global $\GG$-shtukas with fixed legs. Then one obtains analogues/generalizations of Serre-Tate's theorem, as well as of several classical results comparing Shimura varieties and Rapoport-Zink spaces in the arithmetic case. In the last section we define period spaces as suitable subspaces of an affine Grassmannian. Contrary to the arithmetic case they are no longer subspaces of a classical flag variety because we allow also non-minuscule bounds. We define the period map from the analytic space associated with a Rapoport-Zink space to the corresponding period space, and discuss its image and compatibility with the tower of coverings of the generic fiber of the Rapoport-Zink space.

\section{Local $\G$-shtukas}\label{sec2}
\subsection{Generalities}
Let $\BF_q$ be a finite field of characteristic $p$ with $q$ elements, let $\BF$ be a fixed algebraic closure of $\BF_q$, and let $\BF_q\dbl z\dbr$ and $\BF_q\dbl\zeta\dbr$ be the power series rings over $\BF_q$ in the (independent) variables $z$ and\ $\zeta$. As base schemes we will consider the category $\Nilp_{\BF_q\dbl\zeta\dbr}$ consisting of schemes over $\Spec\BF_q\dbl\zeta\dbr$ on which $\zeta$ is locally nilpotent. Let $\G$ be a parahoric group scheme over $\Spec\BF_q\dbl z\dbr$  with connected reductive generic fiber, compare \cite{BT72}, D\'ef. 5.2.6 and \cite{AnhangPR}.

Let $S\in \Nilp_{\BF_q\dbl\zeta\dbr}$ and consider any sheaf of groups $H$ on $S$ for the fpqc-topology. By an $H$-torsor on $S$ we mean a sheaf $\mathcal H$ for the fpqc-topology on $S$ together with a (right) action of the sheaf $H$ such that $\mathcal H$ is isomorphic to $H$ on an fpqc-covering of $S$. 

Let now $L\G$ and $L^+{\G}$ be the loop group and the group of positive loops associated with $\G$, i.e. for an $\BF_q$-algebra $R$ let 
$$(L^+{\G})(R)=\G(R\dbl z\dbr) \quad\text{ and }\quad(L\G)(R)=\G(R\dpl z\dpr).$$

Let $\CG$ be an $L^+{\G}$-torsor on $S$. Via the inclusion $L^+{\G}\subset L\G$ we can associate an $L\G$-torsor $L\CG$ with $\CG$. For any $L\G$-torsor $\CG'$ on $S$ we denote by $\sigma^*\CG'$ the pullback of $\CG'$ under the $q$-Frobenius morphism $\sigma:={\rm Frob}_q: S\rightarrow S$.

\begin{definition}\label{DefLocSht}
A \emph{ local $\G$-shtuka} over some $S\in\Nilp_{\BF_q\dbl\zeta\dbr}$ is a pair $\underline{\CG} = (\CG,\tau_\CG)$ consisting of an $L^+{\G}$-torsor $\CG$ on $S$ and an isomorphism of the associated $L\G$-torsors $\tau_\CG\colon\sigma^\ast L\CG \isoto L\CG$.

A \emph{ quasi-isogeny} $g\colon(\CG',\tau_{\CG'})\to(\CG,\tau_\CG)$ between local $\G$-shtukas over $S$ is an isomorphism $g\colon L\CG'\isoto L\CG$ of the associated $L\G$-torsors with $g\circ\tau_{\CG'}=\tau_\CG\circ\sigma^* g$.
\end{definition}

Local $\G$-shtukas were introduced and studied in joint work of Hartl and the author \cite{HV1}, \cite{HV2} in the case where $\G$ is a constant split reductive group over $\BF_q$. The general case was first considered in work of Arasteh Rad and Hartl \cite{AH_Local}.

\begin{example}\label{exlocsht}
For $\G=\GL_r$, we have the following more classical description.

A local shtuka over $S\in\NilpFq$ (of rank $r$) is a pair $(M,\phi)$ where $M$ is a sheaf of $\mathcal{O}_S\dbl z\dbr$-modules on $S$ which Zariski-locally is free of rank $r$, together with an isomorphism of $\mathcal{O}_S\dpl z\dpr$-modules $$\phi:\sigma^* M\otimes_{\mathcal{O}_S\dbl z\dbr}\mathcal{O}_S\dpl z\dpr\isoto M\otimes_{\mathcal{O}_S\dbl z\dbr}\mathcal{O}_S\dpl z\dpr.$$
Then the category of local $\GL_r$-shtukas over $S$ is equivalent to the category of local shtukas of rank $r$ over $S$, see \cite{HV1}, Lemma 4.2. 

Local shtukas were first introduced by Anderson \cite{Anderson} over a complete discrete valuation ring. Genestier \cite{Genestier} constructed moduli spaces for them in the Drinfeld case and used these to uniformize Drinfeld modular varieties. 
\end{example}

An important invariant of local $\G$-shtukas is their \emph{Newton point}. To define it let $k$ be an algebraically closed field of characteristic $p$ and let $L=k\dpl z\dpr$. Then every local $\G$-shtuka $(\CG, \tau_\CG)$ over $k$ has a trivialization $\CG\cong (L^+{\G})_{k}$. Via this isomorphism the Frobenius map $\tau_{\CG}$ corresponds to an element $b\in L\G(k)=\G(L)$. Changing the trivialization replaces $b$ by a different representative of its $L^+{\G}(k)$-$\sigma$-conjugacy class. In the same way, local $\G$-shtukas isogenous to $(\CG, \tau_\CG)$ correspond to the elements of the $\sigma$-conjugacy class
$$[b]=\{g^{-1}b\sigma(g)\mid g\in \G(L)\}.$$
Let $\Geta$ denote the generic fiber of $\G$. Then the set of $\sigma$-conjugacy classes $B(G)=\{[b]\mid b\in G(L)\}$ for quasi-split $G$ is described by Kottwitz in \cite{Kottwitz1}, \cite{KottwitzIC2} by two invariants. Let $T$ be a maximal torus of $G$. The first invariant is then the Kottwitz point, i.e. the image under the Kottwitz map $\kappa_G: G(L)\rightarrow \pi_1(G)_{\Gamma}$. Here, $\pi_1(G)$ is Borovoi's fundamental group, and $\Gamma$ is the absolute Galois group of $\mathbb F_q\dpl z\dpr$. This invariant is locally constant on $LG$, and can also be computed using the identification $\pi_0(LG)=\pi_1(G)/\Gamma$, compare \cite{Neupert}, 2.2. The second invariant is the Newton point $\nu_b$, an element of $X_*(T)_{\mathbb{Q}}^{\Gamma}$.

Consider now a local $\G$-shtuka $(\CG, \tau_\CG)$ over a scheme $S$. Then we obtain an induced decomposition of $S(k)$ into subsets
 $\mathcal N_{[b]}(k)$ for $[b]\in B(G)$ with $$\mathcal N_{[b]}(k)=\{x\in S(k)\mid (\CG, \tau_\CG)_x \text{ is in the isogeny class of }[b]\}.$$  By \cite{RapoportRichartz}, this induces a decomposition of $S$ into locally closed subschemes, which we call Newton strata. If $S$ is connected, they correspond to the strata obtained by fixing the Newton point only.
 
\subsection{Bounds}

To ensure finiteness properties of moduli spaces of local $\G$-shtukas, or of deformations of local $\G$-shtukas, we bound the singularity of the morphism $\tau_{\mathcal G}$. In the arithmetic context, such bounds are usually assumed to be minuscule and correspond to the choice of the Hodge cocharacter defining a Shimura datum. In our context we define bounds more generally as suitable ind-subschemes of the affine flag variety.

The \emph{affine flag variety $\Flag_{\G}$ associated with $G$} is the fpqc-sheaf associated with the presheaf
$$
S\mapsto L\G(S)/L^+{\G}(S)=\G(\mathcal O_S\dpl z \dpr(S) )/\G(\mathcal O_S\dbl z\dbr(S))
$$ 
on the category of $\BF_q$-schemes. By \cite{PR2}, Theorem 1.4 and \cite{Richarz13} Theorem A, $\Flag_\G$ is represented by an ind-scheme which is ind-projective over $\BF_q$.

Consider further the group scheme $\G\times_{\BF_q\dbl z\dbr}\Spec \BF_q\dpl\zeta\dpr\dbl z-\zeta\dbr$ under the homomorphism 
\begin{equation}\label{eqmorgr}\BF_q\dbl z\dbr\to \BF_q\dpl\zeta\dpr\dbl z-\zeta\dbr,\quad\quad z\mapsto z=\zeta+(z-\zeta).\end{equation} The associated affine Grassmannian $\Gr_{\G}^{B_{\dR}}$ is the sheaf of sets for the fpqc-topology on $\Spec \BF_q\dpl\zeta\dpr$ associated with the presheaf
\begin{equation}\label{EqH_G}
X \mapsto \G(\mathcal{O}_X\dpl z-\zeta\dpr)/\G(\mathcal{O}_X\dbl z-\zeta\dbr).
\end{equation}
\begin{remark}\label{remgrdr}
$\Gr_{\G}^{B_\dR}$ is in the same way as $\Flag_{\G}$ above an ind-projective ind-scheme over $\Spec\BF_q\dpl\zeta\dpr$. Here, note that the homomorphism in \eqref{eqmorgr} induces an inclusion $\BF_q\dpl z\dpr\to \BF_q\dpl\zeta\dpr\dbl z-\zeta\dbr$, so the group $\G\times_{\BF_q\dbl z\dbr}\Spec \BF_q\dpl\zeta\dpr\dbl z-\zeta\dbr$ is reductive, justifying the name affine ``Grassmannian''. The notation $B_{\dR}$ refers to the fact that if $C$ is the completion of an algebraic closure of $\BF_q\dpl\zeta\dpr$, then $C\dpl z-\zeta\dpr$ is the function field analog of Fontaine's $p$-adic period field $B_\dR$.
\end{remark}

We fix an algebraic closure $\overline{\BF_q\dpl\zeta\dpr}$ of $\BF_q\dpl\zeta\dpr$, and consider pairs $(R,\hat Z_R)$, where $R/\BF_q\dbl\zeta\dbr$ is a finite extension of discrete valuation rings contained in $\overline{\BF_q\dpl\zeta\dpr}$, and where $\hat Z_R\subset \widehat\Flag_{\G,R}:=\Flag_{\G} \widehat\times_{\BF_q}\Spf R$ is a closed ind-subscheme that is contained in a bounded (or: projective) subscheme of $\widehat\Flag_{\G,R}$. Two such pairs $(R,\hat Z_R)$ and $(R',\hat Z'_{R'})$ are equivalent if they agree over a finite extension of the rings $R$, $R'$.

A \emph{ bound} is then abstractly defined as an equivalence class $\hat Z:=[(R,\hat Z_R)]$ of such pairs $(R,\hat Z_R)$ satisfying certain properties (for the precise conditions compare \cite{hv3}, Def. 2.1), in particular
\begin{enumerate}
\item $\hat{Z}_R\subset\widehat{\Flag}_{\G,R}$ is a $\zeta$-adic formal scheme over $\Spf R$ that is stable under the left $L^+{\G}$-action.
\item The special fiber $Z_R:=\hat{Z}_R\widehat\times_{\Spf R}\Spec \kappa_R$ is a quasi-compact subscheme of $\Flag_G\widehat\times_{\BF_q}\kappa_R$ where $\kappa_R$ is the residue field of $R$.
\item  Let $\hat{Z}_R^\an$ be the strictly $R[\tfrac{1}{\zeta}]$-analytic space associated with $\hat{Z}_R$. Then the $\hat{Z}_R^\an$ are invariant under the left multiplication of $\G(\cdot\dbl z-\zeta\dbr)$  on $\Gr_{\G}^{\text{{\bf B}}_{\rm dR}}$.
\end{enumerate}
The reflex ring $R_{\hat Z}$ of $[(R,\hat Z_R)]$ is the intersection in $\overline{\BF_q\dpl\zeta\dpr}$ of the fixed field of $\{\gamma\in{\rm Aut}_{\mathbb F_q\dbl\zeta\dbr}(\overline{\mathbb F_q\dpl\zeta\dpr})\mid \gamma(\hat{Z})=\hat{Z}\}$ with all finite extensions of $\BF_q\dbl\zeta\dbr$ over which a representative $\hat{Z}_R$ of the given bound exists.

Let $(\CG, \tau_\CG)$ be a local $\G$-shtuka over some $S\in {\rm Nilp}_{R_{\hat Z}}$ and let $S'$ be an \'etale covering of $S$ over which a trivialization $\CG\cong(L^+{\G})_{S'}$ exists. Then $(\CG, \tau_\CG)$ (or $\tau_{\CG}$) is \emph{bounded by $\hat Z$} if for every such trivialization and for every finite extension $R$ of $\BF_q\dbl\zeta\dbr$ over which a representative $\hat Z_R$ of $\hat Z$ exists, the morphism
\begin{equation*}
 S'\widehat\times_{R_{\hat Z}}\Spf R\rightarrow L\G\widehat\times_{\BF_q}\Spf R\rightarrow \widehat{\Flag}_{\G,R}
\end{equation*}
induced by $\tau_\CG$ factors through $\hat{Z}_R$.

The above definition of bounds follows the strategy to allow as general bounds as possible. On the other hand, classically and by analogy with the arithmetic case, one considers bounds given as Schubert varieties. They are described as follows (compare \cite{AH_Local}, Example 4.12).  

\begin{example}\label{exbound}
Consider the base change $\Geta_L$ of $\Geta$ to $L=\mathbb{F}\dpl z\dpr$. Let $A$ be a maximal split torus in $\Geta_L$ and let $T$ be its centralizer, a maximal torus. Let $N = N(T)$ be the normalizer of $T$ and let $\mathcal{T}^0$ be the identity component of the N\'eron model of $T$ over $\mathcal{O}_L=\mathbb{F}\dbl z\dbr$.

Consider the Iwahori--Weyl group $\widetilde{W}= N(L)/\mathcal{T}^0(\mathcal{O}_L)$ and let $\widetilde W^{\G}=(N(L)\cap \G(\mathcal O_L))/\mathcal{T}^0(\mathcal{O}_L)$.
Then by the Bruhat-Tits decomposition we have a bijection
\begin{equation}\label{eqdecomp}
 \widetilde{W}^{\G}  \backslash \widetilde{W}/ \widetilde{W}^{\G}\rightarrow L^+{\G}(\mathbb{F})\backslash L\G(\mathbb{F})/L^+{\G}(\mathbb{F}).
\end{equation}

Let $\omega$ be in the right hand side, and let $\mathbb{F}_{\omega}$ be a finite extension of $\mathbb{F}_q$ such that $\omega$ has a representative $g_{\omega}\in L\G(\mathbb{F}_{\omega})$. We define the \emph{Schubert variety} $S_\omega$ as the ind-scheme theoretic closure of the $L^+{\G}$-orbit of $g_{\omega}$ in $\Flag_{\G}\widehat\times_{\mathbb{F}_q}\mathbb{F}_{\omega}$. Let $R=\mathbb{F}_{\omega}\dbl\zeta\dbr$ and let $\hat{Z}_{\mathbb{F}_{\omega}\dbl\zeta\dbr}=S_{\omega}\widehat\times_{\mathbb{F}_\omega}\Spf R$. Then the equivalence class of $(R, \hat{Z}_R)$ defines a bound. In this case we also say ``bounded by $\omega$'' instead of ``bounded by $(R, \hat{Z}_R)$''.
\end{example}

\section{Deformations}\label{sec3}

Let $\Flag_{\G}$ denote again the affine flag variety, i.e.~the quotient sheaf $\Flag_{\G}=L\G/L^+{\G}$, an ind-scheme over $\mathbb F_q$ which is of ind-finite type. Let $\widehat \Flag_{\G}$ be the fiber product $\Flag_{\G}\times_{\Spec \mathbb{F}_q}\Spf~\mathbb{F}_q\dpl \zeta\dpr$.

Then generalizing \cite[Prop. 5.1]{HV1} (where the case that $\G$ is split is considered) the ind-scheme $\widehat \Flag_{\G}$ pro-represents the functor $\NilpFq^o \rightarrow (\Sets)$
\begin{align*}
S\mapsto & \{(\mathcal G,\delta)\mid \mathcal G \text{ a $\G$-torsor on $S$,} \\
&\quad\delta:L\mathcal G\rightarrow L\G_S\text{ an isomorphism of the associated $L\G$-torsors}\}/\cong.
\end{align*}
Here $(\mathcal G,\delta)$ and $(\mathcal G',\delta)$ are isomorphic if $\delta^{-1}\circ\delta'$ is an isomorphism $\mathcal G'\to\mathcal G$.

We fix a bound $[(R,\hat Z_R)]$ and a local $\G$-shtuka $\underline{\mathbb{G}}=(\mathbb{G},\phi_{\mathbb G})$ over a field $k\in\NilpFq$ which is bounded by $\hat Z_R$. We consider the functor of bounded deformations of $\underline{\mathbb{G}}$ on the category of Artinian local $k\dbl\zeta\dbr$-algebras with residue field $k$,
\begin{align*}
{\rm Def}_{\underline{\mathbb{G}}, \hat Z_R}:({\rm Art}_{k})\rightarrow& (\Sets)\\
A\mapsto& \{(\underline{\CG},\beta)\mid \underline\CG\text{ a local $\G$-shtuka over $\Spec A$ bounded by $\hat Z_R$}\\
 & \quad \beta:\underline{\mathbb G}\rightarrow\underline\CG\otimes_A k \text{ an isomorphism of local $\G$-shtukas}\}/\cong
\end{align*}
and require isomorphisms to be compatible with the maps $\beta$.

\begin{remark}
In \cite{HV1}, 5 this functor and the deformation space representing it are explicitly described in case that $\G$ is split and $\hat Z_R$ associated with some $\omega\in \widetilde W$ as in Example \ref{exbound}, see below. This is generalized to reductive groups $\G$ in \cite{vwu}. It would be interesting to have these results also in the more general context of parahoric group schemes.
\end{remark}

Assume that $\G$ is a reductive group scheme over $\Spec\BF_q\dbl z\dbr$. Then the left hand side of (\ref{eqdecomp}) is $W_0\backslash\widetilde W/W_0$ where $W_0$ is the finite Weyl group of $\G_{\BF_q}$. This double quotient is isomorphic to $X_*(T)_{\dom}$. Let $\mu\in X_*(T)_{\dom}$ and $(R,\hat Z_R)$ be the corresponding bound as in Example \ref{exbound}. Let $\underline{\mathbb{G}}$ be a local $\G$-shtuka bounded by $\mu$ over a field $k\in\NilpFq$. Assume that there is a trivialization $\alpha:\underline{\mathbb{G}}\rightarrow(\G_{k},b_0\sigma)$ for some $b_0\in L\G(k)$. Using the boundedness, let $x\in Z_R(k)$ be the point defined by $b_0^{-1}$. By \cite{HV1}, Thm. 5.6 and \cite{vwu}, Proposition 2.6 we have 
\begin{thm}
Let $\G$ be reductive, and $\hat Z_R$ as in Example \ref{exbound}. Let $\underline{\mathbb{G}}$ and $\alpha$ be as above. Let $D$ be the complete local ring of $\hat Z_R$ at the point $x$. It is a complete noetherian local ring over $k\dbl\zeta\dbr$.
Then $D$ pro-represents the formal deformation functor ${\rm Def}_{\underline{\mathbb{G}}, \hat Z_R}$.
\end{thm}

To study Newton strata in a family of local $\G$-shtukas it is helpful to also have a corresponding stratification on deformations. However, for infinitesimal deformations such a question does not make sense. Instead we will replace the formal deformation space $\Spf D$ by $\Spec D$, using the following proposition. Here, a linearly topologized $\BF_q\dbl z\dbr$-algebra $R$ is admissible if $R=\underset{\leftarrow}{\lim} R_\alpha$ for a projective system $(R_\alpha,u_{\alpha\beta})$ of discrete rings such that the filtered index-poset has a smallest element $0$, all maps $R\to R_\alpha$ are surjective, and the kernels $I_\alpha:=\ker u_{\alpha,0}\subset R_\alpha$ are nilpotent.

\begin{prop}
Let $R$ be an admissible $\BF_q\dbl z\dbr$-algebra as above with filtered index-poset $\mathbb N_0$. Then the pullback under the natural morphism $\Spf R\to \Spec R$ defines a bijection between local $\G$-shtukas bounded by $\mu$ over $\Spec R$ and over $\Spf R$. 
\end{prop}

\begin{remark} Without a boundedness condition the pullback map is in general only injective. The corresponding result for $p$-divisible groups is shown by Messing and by de Jong. The above proposition is proved for split $\G$ in \cite{HV1}, Proposition 3.16. However, the same proof also shows this assertion for all $\G$ that are reductive over $\Spec\BF_q\dbl z\dbr$.
\end{remark}

\section{Moduli spaces}\label{sec4}

In this section we fix $\G$, an isogeny class of local $\G$-shtukas and a bound $\hat Z$. We then define moduli spaces of local $\G$-shtukas bounded by $\hat Z$ in the same way as Rapoport-Zink define their moduli spaces of $p$-divisible groups in \cite{RZ}.

\begin{definition}
Let $\underline{{\mathbb{G}}}_0$ be a local $\G$-shtuka over $\mathbb{F}$. Let $\hat{Z}=[(R,\hat Z_R)]$ be a bound and denote its reflex ring by $R_{\hat Z}=\kappa\dbl\xi\dbr$. It is a finite extension of $\mathbb F_q\dbl \zeta\dbr$. Set $\breve R_{\hat Z}:=\mathbb{F}\dbl\xi\dbr$, and consider the functor
\begin{align*}
\breve{\mathcal M}: (\Nilp_{\breve R_{\hat Z}})^{\circ} &\rightarrow(\Sets)\\
S&\mapsto  \{(\underline\CG,\bar\delta)\mid\underline{\CG}\text{ a local $\G$-shtuka over $S$ bounded by $\hat{Z}^{-1}$,}\\ 
&\quad\quad\quad\quad\quad \bar{\delta}: \underline{\CG}_{\bar{S}}\rightarrow \underline{\mathbb{G}}_{0,\bar{S}}\text{ a quasi-isogeny}\}/\cong
\end{align*}
Here $\bar{S}:=V_S(\zeta)$ is the zero locus of $\zeta$ in $S$, and two pairs $(\underline\CG,\bar\delta),(\underline\CG',\bar\delta')$ are isomorphic if $\bar{\delta}^{-1}\circ \bar{\delta'}$ lifts to an isomorphism $\underline\CG'\rightarrow\underline\CG$. 
\end{definition}
By \cite{hv3}, 2.2) $\underline{\CG}$ is bounded by $\hat{Z}^{-1}$ if and only if $\tau_{\CG}^{-1}$ is bounded by $\hat Z_R$.

Let ${\rm QIsog}_{\mathbb{F}}(\underline{{\mathbb{G}}}_0)$ be the group of self-quasi-isogenies of $\underline{{\mathbb{G}}}_0$. It naturally acts on the functor $\breve{\mathcal M}$. Since $\mathbb{F}$ has no non-trivial \'etale coverings, we may fix a trivialization $\underline{\mathbb{G}}_0\cong((L^+{\G})_\mathbb{F},b\sigma^*)$ where $b\in L\G(\mathbb{F})$ represents the Frobenius morphism. Via such a trivialization, ${\rm QIsog}_\mathbb{F}(\underline{\mathbb{G}}_0)$ is identified with 
\begin{equation}\label{eqdefj}
J_b(\mathbb{F}_q\dpl z\dpr):=\{g \in \Geta(\mathbb{F}\dpl z\dpr)\mid g^{-1}b\sigma(g)=b\}. 
\end{equation}
This is the set of $\mathbb{F}_q\dpl z\dpr$-valued points of an algebraic group $J_b$ over $\mathbb{F}_q\dpl z\dpr$.

\begin{thm}[\cite{AH_Local}, Thm.~4.18 and Cor.~4.26, \cite{hv3}, Rem. 3.5]\label{thmrzrep}
The functor $\breve{\mathcal M}$ is ind-representable by a formal scheme over $\Spf \breve R_{\hat Z}$ which is locally formally of finite type and separated. It is an ind-closed ind-subscheme of $\Flag_{\G}\widehat\times_{\mathbb{F}_q}\Spf\breve R_{\hat Z}$.
\end{thm}

By its analogy with moduli spaces of $p$-divisible groups, the formal scheme representing $\breve{\mathcal M}$ is called a Rapoport-Zink space for bounded local $\G$-shtukas.

Let $E$ be the quotient field of $R_{\hat Z}$, and $\breve E\cong \mathbb{F}\dpl \xi\dpr$ the completion of its maximal unramified extension. We write ${\breve{\mathcal M}}^{\an}$ for the strictly $\breve E$-analytic space associated with $\breve {\mathcal M}$. 

Next we explain the construction of a tower of coverings of ${\breve{\mathcal M}}^{\an}$. Roughly spoken, it is obtained by trivializing the (dual) Tate module of the universal local $\G$-shtuka over $\breve{\mathcal M}^{\an}$. 

\begin{definition}
Let $S$ be an $\BF\dpl\zeta\dpr$-scheme or a strictly $\BF\dpl\zeta\dpr$-analytic space. Then an \emph{\'etale local $\G$-shtuka} over $S$ is a pair $\underline{\CG}=(\CG,\tau_{\CG})$ consisting of an $L^+{\G}$-torsor $\CG$ on $S$ and an isomorphism $\tau_\CG:\sigma^*\CG\rightarrow\CG$ of $L^+{\G}$-torsors.
\end{definition}

\begin{remark}
Let $S=\Spf B$ be an affinoid admissible formal $\breve{R}_{\hat Z}$-scheme and $S^{\an}$ the associated strictly $\breve E$-analytic space. Consider a trivialized local $\G$-shtuka $((L^+{\G})_{S},A\sigma^*)$ over $S$. Then $A\in \G(B\dbl z\dbr[\tfrac{1}{z-\zeta}])$. Note that $$(z-\zeta)^{-1}=-\sum_{i=0}^\infty\zeta^{-i-1}z^i\in\mathcal{O}_{S^{\an}}(S^{\an})\dbl z\dbr$$ implies $B\dbl z\dbr[\tfrac{1}{z-\zeta}]\subset\mathcal{O}_{S^{\an}}(S^{\an})\dbl z\dbr$. Therefore $((L^+{\G})_{S},A\sigma^*)$ induces an \'etale local $\G$-shtuka $((L^+{\G})_{S^{\an}},A\sigma^*)$ over $S^{\an}$.

To obtain similarly a  universal family of \'etale local $\G$-shtukas over ${\breve{\mathcal M}}^{\an}$, we cover $\breve{\mathcal M}$ by affinoid admissible formal $\breve{R}_{\hat Z}$-schemes. For each of them one chooses a finite \'etale covering trivializing the local $\G$-shtuka, and applies the above construction. Descending the \'etale local $\G$-shtuka then yields the desired universal family, compare \cite{hv3}, 6.
\end{remark}

There are two approaches to construct the (dual) Tate module of an \'etale local $\G$-shtuka $\underline{\CG}$ over a connected strictly $\breve E$-analytic space $X$. To define it directly, following Neupert \cite{Neupert}, 2.6, consider for each $n\in\mathbb{N}$ the $L^+{\G}/\G_n$-torsor associated with $\CG$ where $\G_n$ is the kernel of the projection $\G(\BF_q\dbl z\dbr)\rightarrow \G(\BF_q\dbl z\dbr/(z^n))$. The isomorphisms $\tau_\CG:\sigma^*\CG\rightarrow\CG$ and $\sigma^{*}:\CG\rightarrow\sigma^*\CG$ then induce corresponding maps of $L^+{\G}/\G_n$-torsors. The invariants of $\tau\circ \sigma^*$ form a $\G(\BF_q\dbl z\dbr/(z^n))$-torsor which is trivialized by a finite \'etale covering of $X$. One can then define the Tate module of $\underline\CG$ as the inverse limit over $n$ of these torsors. 

Alternatively, one can define it as a tensor functor, following \cite{AH_Local}, or \cite{hv3}, 7: Fix a geometric base point $\bar x$ of $X$. We consider representations $\rho:\G\rightarrow \GL_r$ in $\Rep_{\mathbb{F}_q\dbl z\dbr}(\G)$. Let $\underline M=(M,\tau_M)$ be the \'etale local shtuka of rank $r$ associated with $\rho_*\underline{\CG}$ as in Example \ref{exlocsht}. Let $\underline M_{\bar x}$ denote its fiber over $\bar x$. The (dual) Tate module of $\underline{\mathcal G}$ with respect to $\rho$ is the (dual) Tate module of $\underline M_{\bar x}$, $$\check T\underline M_{\bar x}=\{m\in \underline M_{\bar x} \mid \tau_M(\sigma^*m)=m\},$$
a free $\mathbb F_q\dbl z\dbr $-module of rank $r$. It carries a continuous monodromy action of $\pi_1^\et(X,\bar x)$, which also factors through $\pi_1^\alg(X,\bar x)$. We obtain the \emph{dual Tate module} of $\underline{\mathcal G}$ as a tensor functor 
$$\check T_{\underline{\CG},\bar x}:\Rep_{\BF_q\dbl z\dbr}\G\rightarrow\Rep^{\rm cont}_{\BF_q\dbl z\dbr }(\pi_1^{\rm alg}(X,\bar x)).$$
Similarly, using rational representations in $\Rep_{\BF_q\dpl z\dpr}\G$ one can define the \emph{rational dual Tate module}
$$\check V_{\underline{\CG},\bar x}:\Rep_{\BF_q\dpl z\dpr}\G\rightarrow \Rep^{\rm cont}_{\BF_q\dpl z\dpr }(\pi_1^{\rm et}(X,\bar x)).$$
The two constructions of dual Tate modules are compatible in the sense that the tensor functor associated with the above torsor coincides with the tensor functor $\check T_{\underline{\CG},\bar x}:\Rep_{\BF_q\dbl z\dbr}\G\rightarrow(\mathbb{F}_q\dbl z\dbr\text{-Loc})_X$ with values in the local systems of $\mathbb{F}_q\dbl z\dbr$-lattices on $X$ that can be associated with $\check T_{\underline{\CG},\bar x}$ as in the second construction (compare \cite{hv3}, Prop. 5.3). 
                                                                                                                                    
We can now proceed to define level structures and the tower of coverings of $\breve{\mathcal M}^{\an}$. Let $\FMod_A$ denote the category of finite locally free $A$-modules. We consider the forgetful functors 
\begin{equation}\label{eqw0}
\omega^{\circ}_A:\Rep_A \G\rightarrow\FMod_A
\end{equation}
 and 
\begin{equation}\label{eqforget}forget:\Rep^{\rm cont}_{A}(\pi_1^{\rm et}(X,\bar x))\rightarrow\FMod_A.
\end{equation}
For an \'etale local $\G$-shtuka $\underline{\CG}$ over $X$ the two sets 
\begin{align*}
{\rm Triv}_{\underline\CG,\bar x}(\BF_q\dbl z\dbr) & =  {\rm Isom}^{\otimes}(\omega^{\circ}_{\BF_q\dbl z\dbr },{\it forget}\circ\check T_{\underline{\CG},\bar x})(\BF_q\dbl z\dbr)\\
{\rm Triv}_{\underline\CG,\bar x}(\BF_q\dpl z\dpr) & ={ \rm Isom}^{\otimes}(\omega^{\circ}_{\BF_q\dpl z\dpr },{\it forget}\circ\check V_{\underline{\CG},\bar x})(\BF_q\dpl z\dpr)
\end{align*}
are non-empty and carry natural actions of $\G(\BF_q\dbl z\dbr )\times \pi_1^{\rm alg}(X,\bar x)$ and $\Geta(\mathbb F_q\dpl z\dpr)\times \pi_1^{\rm et}(X,\bar x)$, respectively. Here, the first factor acts via $\omega^{\circ}$ and the second via the action on the Tate module. Note that non-emptiness needs our assumption that $\G$ has connected fibers, and thus slightly differs from the setting of Rapoport and Zink. 

In the same way as in the arithmetic case we define coverings of $\breve{ \mathcal M}$.
\begin{definition}\label{DefLevel}
Let $\underline{\CG}$ be an \'etale local $\G$-shtuka over a connected $\breve E$-analytic space $X$, and let $K$ be an open compact subgroup of $\G(\BF_q\dbl z\dbr)$ resp. $G(\BF_q\dpl z\dpr)$. An \emph{integral resp. rational $K$-level structure} on $\underline{\CG}$ is a $\pi_1^{\rm alg}(X,\bar x)$- resp.~$\pi_1^{\rm et}(X,\bar x)$-invariant $K$-orbit in ${\rm Triv}_{\underline{\CG},\bar x}(\BF_q\dbl z\dbr)$.

For an open subgroup $K\subset \G(\BF_q\dbl z\dbr)$ let $X^K$ be the functor on the category of $\breve E$-analytic spaces over $X$ parametrizing integral $K$-level structures on the  local $\G$-shtuka $\underline\CG$ over $X$. 

Let $K\subset \G(\BF_q\dpl z\dpr)$ be compact open. Let $K'\subset K$ be a normal subgroup of finite index with $K'\subset \G(\BF_q\dbl z\dbr)$. Then $gK'\in K/K'$ acts on $X^{K'}$ by sending the triple $(\underline\CG,\delta,\eta K')$ over $X^{K'}$ to the triple $(\underline\CG,\delta,\eta gK')$.

We define ${\breve{\mathcal M}}^K$ as the $\breve E$-analytic space which is the quotient of ${\breve{\mathcal M}}^{K'}:={\breve{\mathcal M}}^{{\rm an},K'}$ by the finite group $K/K'$. It is independent of the choice of $K'$. 
In particular, $\breve{\mathcal M}^{K_0}=(\breve{\mathcal M})^{\rm an}$ for $K_0=\G(\BF_q\dbl z\dbr)$. 
\end{definition}

By \cite{hv3}, Cor. 7.13 one can also describe ${\breve{\mathcal M}}^K$ directly as a parameter space of local $\G$-shtukas with level structure.


Furthermore, if $K'\subseteq K$ is as above, the action of $gK'\in K/K'$ on ${\breve{\mathcal M}}^{K'}$ can be described in this interpretation as a Hecke correspondence.

\section{The geometry of the special fiber}\label{sec6}

We fix again a bound $(R,\hat Z_R)$ and let $\kappa$ be the residue field of $R$. Let $Z^{-1}$ be the special fiber of the inverse bound $\hat Z^{-1}$ over $\kappa$. The affine Deligne-Lusztig variety associated with an element $b\in L\G(\mathbb{F})$ and $Z^{-1}$ is the reduced closed ind-subscheme $X_{Z^{-1}}(b)\subset\Flag_{\G}$ whose $k$-valued points (for any field extension $k$ of $\mathbb{F}$) are given by
\begin{equation}\label{defadlv}
X_{Z^{-1}}(b)(k)=\{ g\in \Flag_{\G}(k)\mid g^{-1}b\sigma(g) \in Z^{-1}(k)\}.
\end{equation}
Our conditions on the bound imply that $Z^{-1}$ is a left-$L^+\G$-invariant subscheme of $\Flag_{\G}$. The $L^+\G$-orbits in $\Flag_{\G}$ correspond bijectively to the elements of $\widetilde{W}^{\G}  \backslash \widetilde{W}/ \widetilde{W}^{\G}$, compare Example \ref{exbound}. For $x\in \widetilde W$ let 
\begin{equation}
X_{x}(b)(k)=\{ g\in \Flag_{\G}(k)\mid g^{-1}b\sigma(g) \in L^+{\G} xL^+{\G}\}.
\end{equation}
Thus by left invariance every $X_{Z^{-1}}(b)(k)$ is a union of affine Deligne-Lusztig varieties of the form $X_x(b)$. By the boundedness condition on $\hat Z_R$, this union is finite. In the arithmetic context, such unions of affine Deligne-Lusztig varieties have been studied for example in \cite{HeKR}.

Affine Deligne-Lusztig varieties are the underlying reduced subschemes of Rapoport-Zink spaces:
\begin{thm}[\cite{HV1}, 6, \cite{AH_Local}, Thm.~4.18, Cor.~4.26]
The underlying reduced subscheme of the moduli space $\breve{\mathcal M}$ associated with $\hat Z_R$ and $b$ as in Theorem \ref{thmrzrep} equals $X_{Z^{-1}}(b)$. It is a scheme locally of finite type and separated over $\mathbb{F}$, all of whose irreducible components are projective.
\end{thm}

Thus to describe the special fiber of the Rapoport-Zink spaces $\breve{\mathcal M}$ it is enough to study the affine Deligne-Lusztig varieties $X_x(b)$. For an overview of the current state of the art of this field compare \cite{overview_newton} and X. He's talk at this ICM \cite{HeICM}. We will here only discuss one class of results that were obtained as an application of the relation to local $\G$-shtukas, and that were one of the initial goals of this theory. For the following result it is essential to assume that $\G$ is reductive over $\mathbb{F}_q\dbl z\dbr$. For more general groups there are counterexamples to all assertions in the theorem. Recall that for reductive groups $\widetilde{W}^{\G}  \backslash \widetilde{W}/ \widetilde{W}^{\G}\cong X_*(T)_{\dom}$. Let ${\rm def}(b)={\rm rk} ~G - {\rm rk}_{\mathbb{F}_q\dpl z\dpr}J_b$ where $J_b$ is as in \eqref{eqdefj}. Note that $X_*(T)_{\dom}$ is partially ordered by the Bruhat ordering. Via the induced ordering on the Newton points (requiring equality of Kottwitz points) also $B(G)$ inherits an ordering. Let $l$ denote the length of a maximal chain between two elements in the partially ordered set $B(G)$, and let $\rho$ be the half-sum of the positive roots of $\G$.

\begin{thm}[\cite{grothconj}, \cite{irrmin}, 3.2]\label{thm2}
Let $\mu_1\preceq\mu_2\in X_*(T)$ be dominant coweights and $[b]\in B(G)$. 
\begin{enumerate}
\item Let $$S_{\mu_1,\mu_2}=\bigcup_{\mu_1\preceq\mu'\preceq \mu_2} L^+\G\mu'(z)L^+\G.$$ Assume that $\kappa_G(b)=\mu_2$ in $\pi_1(\G)_{\Gamma}$. Then it is pure of codimension $$\l([b],[\mu_2])=\langle \rho,\mu_2-\nu_b\rangle+\frac{1}{2}{\rm def}(b)$$ in $S_{\mu_1,\mu_2}$. The closure of $\mathcal{N}_{[b]}$ in $S_{\mu_1,\mu_2}$ is the union of all $\mathcal{N}_{[b']}$ for $[b']$ with $\kappa_G(b')=\mu_2$ in $\pi_1(\G)_{\Gamma}$ and $\nu_{b'}\preceq\nu_{b}$.
\item $X_{\mu_2}(b)$ and $X_{\preceq\mu_2}(b)=\bigcup_{\mu\preceq\mu_2}X_{\mu}(b)$ are equidimensional of dimension $$\dim X_{\mu_2}(b)=\dim~X_{\preceq\mu_2}(b)=\langle\rho,\mu_2-\nu_b\rangle-\frac{1}{2}{\rm def}(b).$$
\end{enumerate}
\end{thm}
Here, the first assertion is most useful for $\mu_1=\mu_2$, or for $\mu_1$ the unique minuscule coweight with $\mu_1\preceq\mu_2.$

Note that one needs to define the notions of codimension and of the closure of the infinite-dimensional schemes $\mathcal{N}_b$. Both of these definitions use that there is an open subgroup $H$ of $L\G^+$ such that the Newton point of an element of $S_{\mu_1,\mu_2}$ only depends on its $H$-coset, an element of the finite-dimensional scheme $S_{\mu_1,\mu_2}/H$, compare \cite{grothconj}, 4.3. These (co)dimensions are directly linked to the dimensions of Newton strata in the universal deformation of a local $\G$-shtuka as in Section \ref{sec3}. In fact, Theorem \ref{thm2} is equivalent to a completely analogous result on dimensions and closure relations in that context, and this equivalence is also used in the proof of Theorem \ref{thm2}.

The proofs of the two parts of the theorem are closely linked: One first shows the upper bound on $\dim X_{\preceq\mu_2}(b)$ and how to compute $\dim \mathcal{N}_{[b]}$ from $\dim X_{\preceq\mu_2}(b)$. Then a purity result for the Newton stratification proves the lower bound and equidimensionality. For the proof compare \cite{irrmin}, 3.2 explaining how the proof for split $\G$ in \cite{grothconj} can be generalized. The theorem also inspired a corresponding theory for Newton strata in Shimura varieties of Hodge type by Hamacher.

\section{Comparison to global $\GG$-shtukas}\label{sec5}
\subsection{Global $\GG$-shtukas }
Moduli spaces of bounded global $\GG$-shtukas are the function field analogue of Shimura varieties. Using these moduli spaces V.~Lafforgue \cite{VLaf} showed one direction of the global Langlands correspondence for all groups $\GG$ in the function field case. In this section we describe the relation between moduli spaces of bounded local $\G$-shtukas and global $\GG$-shtukas.

To define them, let $C$ be a smooth projective geometrically irreducible curve over $\BF_q$, and let $\GG$ be a parahoric group scheme over $C$, i.e.~a smooth affine group scheme with connected fibers whose generic fiber is reductive over $\BF_q(C)$ and such that for every point $v$ of $C$ for which the fiber above $v$ is not reductive, the group scheme $\GG_v$ is a parahoric group scheme over $A_{v}$, as defined by Bruhat and Tits. Here $A_{v}$ denotes the completion of $k[C]$ at $v$.

\begin{definition}
A \emph{global shtuka} with $n$ legs over a scheme $S$ is a tuple $(\mathscr{G}, s_1,\dotsc,s_n,\varphi)$ where
\begin{enumerate}
\item $\mathscr G\in \mathcal{H}^1(C,\GG)(S)$ is a $\mathfrak G$-torsor over $C\times_{\mathbb F_q}S$
\item the $s_i\in C(S)$ are pairwise disjoint $S$-valued points, called the legs 
\item $\varphi$ is an isomorphism $$\varphi: \sigma^*\mathscr G|_{C\times_{\mathbb F_q}S\setminus\{\bigcup_i \Gamma_{s_i}\}}\rightarrow\mathscr G|_{C\times_{\mathbb F_q}S\setminus\{\bigcup_i \Gamma_{s_i}\}},$$ called the Frobenius isomorphism.
\end{enumerate} 

The stack $\nabla_n\mathcal{H}^1(C,\GG)$ over $\Spec \BF_q$ is the \emph{moduli stack of global $\mathfrak G$-shtukas with $n$ legs}.
\end{definition}

There is a canonical morphism $\nabla_n\mathcal{H}^1(C,\GG)\rightarrow C^n\setminus \Delta$ mapping a global $\GG$-shtuka to its legs.

We fix distinct places $c_i\in C(\mathbb{F})$ for $1\leq i\leq n$, and write $c=(c_1,\dotsc, c_n)$. Let $A_{c}$ be the completion of the local ring $\mathcal O_{C^n,c}$, and let $\mathbb F_{c}$ be the residue field. Then $A_c\cong \mathbb F_c\dbl \zeta_1,\dotsc, \zeta_n\dbr$.

Writing $A_{c_i}\cong \mathbb{F}_{c_i}\dbl z\dbr$, let $\G_{c_i}=\GG\times_C\Spec A_{c_i}$, a parahoric group scheme over $\Spec \mathbb{F}_{c_i}\dbl z\dbr$, and $\G_i:=\Res_{\mathbb{F}_{c_i}\dbl z\dbr/\mathbb{F}_q\dbl z\dbr}\G_{c_i}$.

Let $$\nabla_n\mathcal H^1(C,\GG)^{c}=\nabla_n\mathcal H^1(C,\GG)\widehat\times_{C^n}\Spf A_{c}$$
be the formal completion of the stack $\nabla_n\mathcal H^1(C,\GG)$ along $c$. It parametrizes global $\GG$-shtukas with $n$ legs in the formal neighborhoods of the $c_i$. We want to define a global-local functor associating with such a global $\GG$-shtuka local shtukas at each of these places, compare \cite{AH_Local}, or \cite{Neupert}, 3.2. Let $(\mathscr G, (s_i),\varphi)\in \nabla_n\mathcal H^1(C,\GG)^{c}(S)$. In other words, $s_i: S\rightarrow C$ factors through $\Spf A_{c_i}$. We consider each place $c_i$ separately. Let $S_{c_i}=\Spec \mathbb{F}_{c_i}\times_{\Spec \mathbb F_q}S$. The base change of $\mathscr G$ to $\Spf A_{c_i}\widehat\times_{\mathbb F_q}S=S_{c_i}\dbl z\dbr$ defines an element of $\mathcal H^1(\mathbb{F}_q,\G_{c_i})(S_{c_i}\dbl z\dbr)$. We now use the following translation, cf. \cite{HV1}, Prop. 2.2a, \cite{Neupert}, Prop. 3.2.4, or \cite{AH_Local}, 2.4.

\begin{prop}
Let $S$ be a scheme over $\mathbb F_q$. There is a natural equivalence of categories $$\mathcal H^1(\mathbb{F}_q,\G)(S\dbl z\dbr)\rightarrow \mathcal H^1(\mathbb{F}_q,L^+{\G})(S).$$ 
\end{prop}
We thus obtain an element of $$\mathcal H^1(\mathbb{F}_q,L^+{\G_{c_i}})(S_{c_i} )=\mathcal H^1(\mathbb{F}_q,L^+{\G_{i}})(S ),$$ or an $L^+{\G}_i$-torsor $\mathscr G_i$ over $S$.

Also, the Frobenius morphism $\varphi$ induces local Frobenius morphisms $\varphi_i$ on the $\G_i$-torsors for each $i$. Altogether we obtain a functor 
$$
\mathfrak{L}=(\mathfrak{L}_1,\dotsc, \mathfrak{L}_n):\nabla_n\mathcal H^1(C,\GG)^c\rightarrow \prod_{i}Sht_{\G_i},\quad\quad
(\mathscr G,(s_i),\varphi)\mapsto ((\mathscr G_i, \varphi_i)) 
$$
called the global-local functor.

In the same way, one can associate with every global $\GG$-shtuka $\underline {\mathscr G}\in\nabla_n\mathcal H^1(C,\GG)^{c}(S)$ for some scheme $S$ and every fixed place $s\in C\setminus \{c_i\}$ an \'etale local $\G_s$-shtuka $\mathfrak{L}_s(\underline {\mathscr G})=\underline{\mathscr G}_s$. One can then define level structures on global $\GG$-shtukas away from the legs using the corresponding notion for local $\G$-shtukas. For a detailed discussion compare \cite{Neupert}, 3.4.

\begin{definition}
Let $\underline{\mathscr G}=(\mathscr G,(s_i),\varphi)$ be a global $\GG$-shtuka over $S$, let $D_0\subset C$ be a finite reduced subscheme with $s_i\in (C\setminus D_0)(S)$ for every $i$. For every $v\in D_0$ fix an open subgroup $U_v\subseteq L^+G_v(\mathbb F_q)$ and let $U=\prod U_v$. Then an \emph{integral $U$-level structure} on $\underline{\mathscr G}$ consists of an integral $U_v$-level structure of $\mathfrak{L}_v(\underline {\mathscr G})$ for every $v\in D_0$.

We denote by $\nabla\mathcal H^1_U(C,\GG)^c$ the \emph{stack of global $\GG$-shtukas} in  $\nabla\mathcal H^1(C,\GG)^c$ with $U$-level structure.
\end{definition}

Rather ad hoc boundedness conditions are defined as follows (compare \cite{AH_Unif}).
\begin{definition}
Let  $\hat{Z}_{c}=(\hat{Z}_i)_i$ be a tuple of bounds $\hat{Z}_i\subset\widehat{\Flag}_{\G_{i}}$. Then $(\mathscr G, (s_i),\varphi)\in \nabla_n\mathcal H^1(C,\GG)^{c}(S)$ is \emph{bounded by $\hat{Z}_{c}$} if for every $i$ the associated local $\G_i$-shtuka $(\mathscr G_i, \varphi_i)$ is bounded by $\hat{Z}_i$. We denote by $\nabla_n^{\hat Z_c}\mathcal H^1(C,\GG)^{c}$ the substack of $\nabla_n\mathcal H^1(C,\GG)^{c}$ of global $\GG$-shtukas bounded by $\hat{Z}_{c}$, and analogously for $\nabla_n^{\hat Z_c}\mathcal H^1_U(C,\GG)^{c}$. 
\end{definition}

There are two other definitions of boundedness: A (seemingly more natural) one by Varshavsky \cite{Varshavsky}, Def. 2.4b for which, however, it is not clear how to establish a compatibility between bounded global and local shtukas. In \cite{Neupert}, 3.3 Neupert gives a global definition of boundedness, generalizing that of Varshavsky and he shows that it coincides with the one presented here.
 
\subsection{The Serre-Tate theorem}

Classically, the Serre-Tate theorem gives an equivalence of categories between deformations of an abelian variety over a scheme on which $p$ is locally nilpotent, and its $p$-divisible group. By work of Arasteh Rad and Hartl, the analog in our situation also holds, taking into account that shtukas have several legs, and requiring them to be in the formal neighborhoods of fixed places of $C$. 

Let $S$ be in $\Nilp_{A_{c}}$ and let $j:  \overline{S}\rightarrow S$ be a closed subscheme defined by a locally nilpotent sheaf of ideals. Let $c_1,\dotsc, c_n\in C$ and let $(\mathscr G, (s_i), \varphi)\in \nabla_n\mathcal H^1(C,\GG)^{c}(\overline S)$. Let $((\mathscr G_i, \varphi_i))=\mathfrak{L}(\mathscr G, (s_i), \varphi)$. 

Let ${\rm Def}_S((\mathscr G, (s_i), \varphi))$ be the category of lifts of $(\mathscr G, (s_i), \varphi)$ to $S$, i.e. of pairs $(\underline{\mathscr H},  \alpha)$ where $\underline{\mathscr H}\in \nabla_n\mathcal H^1(C,\GG)^{c}(S)$ and where $\alpha$ is an isomorphism between $(\mathscr G, (s_i), \varphi)$ and the base change of $\underline{\mathscr H}$ to $\overline S$. Morphisms of lifts are isomorphisms compatible with the morphisms $\alpha$. For the local $\G_i$-shtukas we define analogously a category of lifts ${\rm Def}_S((\mathscr G_i, \varphi_i))$.

\begin{thm}[\cite{AH_Local}, Thm. 5.10]\label{Serre-Tate} 
Let $(c_i)_i\in C^n$, let $(\mathscr G, (s_i), \varphi)\in \nabla_n\mathcal H^1(C,\GG)^{c}(\overline S)$ and let $((\mathscr G_i, \varphi_i))=\mathfrak{L}(\mathscr G, (s_i), \varphi)$. Then the global-local functor induces an equivalence of categories
${\rm Def}_S((\mathscr G, (s_i), \varphi))\rightarrow \prod_i {\rm Def}_S((\mathscr G_i, \varphi_i)).$ 
\end{thm}

\subsection{Foliations}

We begin by considering the uniformization morphism for Newton strata by Rapoport-Zink spaces, following \cite{Neupert}, 5 and \cite{AH_Local}, 5, and paralleling the uniformization result of Rapoport and Zink for $p$-divisible groups, \cite{RZ}. In this subsection we assume that $\GG$ is the base change to $C$ of a reductive group $G$ over $\mathbb{F}_q$. It would be interesting to have a generalization of this result to the more general context considered before.

Let $T$ denote a maximal torus of $G$, and fix a Borel subgroup $B$ containing $T$. Let $C$ be as above and fix characteristic places $c_1,\dotsc, c_n\in C$. Let $\mu_i\in X_*(T)$ be dominant cocharacters defined over a finite extension $E$ of $\mathbb{F}_q$, and $\mu=(\mu_1,\dotsc, \mu_n)\in X_*(T)^n$. Fix decent local $\G_{i}$-shtukas $\underline{\mathbb{G}_i}=(L^+{\G}_{i E},b_i\sigma^*)$ over $\Spec E$, in particular $b_i\in L\G_i(E)$. Let $\mathcal M_{b_{i}}^{\preceq \mu_i}$ be the Rapoport-Zink space associated with $[b_i]$ and the bound given by $\mu_i$.

\begin{thm}[\cite{Neupert}, Theorem 5.1.18, \cite{AH_Unif}, Theorem 7.4]\label{thmunif}
 Let $S$ be a DM-stack over $\Spf E\dbl \zeta_1,\dotsc,\zeta_n\dbr$ such that $(\zeta_1,\dotsc,\zeta_n)$ is locally nilpotent on $S$. Let $(\mathscr G_0,\phi_0,\psi_0)\in \nabla_n\mathcal H_U^1(C,\GG)^c(S)$ for some congruence subgroup $U$. For each place $c_i$ assume that there is, and fix, an isomorphism $\mathfrak{L}_{c_i}(\mathscr G_0,\phi_0)\cong \underline{\mathbb{G}_i}_S$. Then there is a morphism of formal DM-stacks over $\Spf E\dbl \zeta_1,\dotsc,\zeta_n\dbr$  
$$S\times_{\Spf E\dbl \zeta_1,\dotsc,\zeta_n\dbr}\prod_{i}\mathcal M_{b_{i}}^{\preceq \mu_i}\rightarrow \nabla_n^{\mu}\mathcal{H}^1_U(C,\GG)^c$$
which is ind-proper and formally \'etale.
\end{thm}
This theorem has also generalizations to coverings of the moduli spaces associated with compatible level structures, and is then equivariant with respect to the action of $\GG(\mathbb{A}^{c_i}_{\mathbb Q})$ by Hecke correspondences. 

We consider for the universal global $\GG$-shtuka $\underline{\mathscr G}$ over $\nabla_n^{\mu}\mathcal{H}^1_U(C,\GG)^c$  the strata
$$\mathcal{N}^{([b_i])_i}_{\underline{\mathscr G}}=\bigcap_i \mathcal{N}_{[b_i],\mathfrak{L}_i(\underline{\mathscr G})}\subseteq \nabla_n^{\mu}\mathcal{H}^1_U(C,\GG)^c.$$ 
Here, the $\mathcal{N}_{[b_i],\mathfrak{L}_i(\underline{\mathscr G})}$ are the Newton strata in $\nabla_n^{\mu}\mathcal{H}^1_U(C,\GG)^c$ associated with the local $\G_i$-shtuka $\mathfrak{L}_i(\underline{\mathscr G})$.

Then the uniformization morphism of Theorem \ref{thmunif} maps each geometric point of the left hand side to $\mathcal{N}^{([b_i])_i}_{\underline{\mathscr G}}$ where the $[b_i]$ are the classes of the elements used to define the Rapoport-Zink space. Refining this description we now consider a foliation structure on these Newton strata, following Neupert \cite{Neupert}, 5. 

In the arithmetic case, Oort and later Mantovan  \cite{MantThesis} defined a foliation structure on Newton strata on Shimura varieties, by Rapoport-Zink spaces and so-called Igusa varieties. The latter are defined as covers of a central leaf, i.e. the locus in the Shimura variety where the $p$-divisible group is in a very particular isomorphism class. In our situation, the definition is slightly more involved. One needs to pass to the perfection of all involved moduli spaces to define the analog of Igusa varieties, and then to construct the foliation morphism. 

Recall that a Newton stratum only depends on $L\G_i$-$\sigma$-conjugacy classes $[b_i]$, and not on individual representatives. A fundamental alcove in $[b_i]$ is an element $x_{b_i}$ of $\widetilde W$ such that all (or one) representative $b_{i,0}$ in $N(L)$ is contained in $[b_i]$, and such that the length $\ell(x_{b_i})$ is minimal with that property, compare \cite{trunc1}, Theorem 6.5, or \cite{Nie}. One can then show that this length is equal to $\langle 2\rho_G,\nu_{b_i}\rangle$ where $\rho_G$ is the half-sum of the positive roots of $G$ and where $\nu_{b_i}$ is the dominant Newton point of $[b_i]$. Fundamental alcoves are a group-theoretic generalization of the minimal $p$-divisible groups studied by Oort. For each $[b_i]$ we fix such a representative $b_{i,0}$. The central leaf $$\mathcal C^{(b_{i,0})}_U\subseteq \mathcal{N}^{([b_i])_i}_{\underline{\mathscr G}}$$ is then defined as the locus where the associated local $\G_i$-shtukas are isomorphic to $((L^+{\G_i})_k,b_{i,0}\sigma^*)$. It depends on the choice of the fundamental alcove, but not on that of the representative in $N(L)$.

We now define formal versions of these subschemes of $\nabla_n^{\mu} \mathcal H^1_U(C,G)$.

\begin{definition}
Recall that the formal completion of $C^n\setminus (c_1,\dotsc, c_n)$ is isomorphic to $\Spf A_c=\Spf ~\mathbb F_c\dbl \zeta_1,\dotsc, \zeta_n\dbr$.

Let $\mathfrak N_U^{([b_i])_i}$ denote the formal completion of $\nabla_n^{\mu} \mathcal H^1_U(C,G)^c$ along the Newton stratum $\mathcal{N}^{([b_i])_i}_{\underline{\mathscr G}}$. This is not $(\zeta_1,\dotsc, \zeta_n)$-adic.

Let $\mathfrak C_U^{(b_{i,0})}$ be the locus in $\nabla_n^{\mu} \mathcal H^1_U(C,G)^c$ such that after an fpqc-covering each local $\G_{i}$-shtuka associated with the universal global $\GG$-shtuka is isomorphic to $(L^+\G_{i},b_{i,0}\sigma^*)$. One can then show that this is represented by a $(\zeta_1,\dotsc, \zeta_n)$-adic formal scheme over $\Spf A_c$, called the \emph{central leaf}. (\cite{Neupert}, Prop. 6.1.7)

By $X^{\sharp}$ we denote the perfection of a formal scheme $X$. It is again a formal scheme.

For $d\in\mathbb{N}$ we consider the subgroups $K_d=\{g\in L^+\G_i\mid g\equiv 1\pmod{z^{d+1}}\}$ and $$I_{d}(b_{i,0})=\bigcap_{N\geq 0}\phi^N(K_d)$$ where $\phi(g)=b_{i,0}^{-1}\sigma^{-1}(g)b_{i,0}$.
\end{definition}

The following theorem shows existence of the formal Igusa varieties.

\begin{thm}[\cite{Neupert}, Cor. 6.1.9]
For all tuples $(d_i)_i$, define $\mathfrak{Ig}_U^{(d_i)\sharp}(T)$ as the set of pairs consisting of an element of  $\mathfrak{C}_U^{(b_i)\sharp}(T)$ together with $I_{d_i}$-truncated isomorphisms between the associated local $G_i$-shtukas and $(L^+G_i,b_{i,0}\sigma^*)$ for all $i$. Then this defines a sheaf $\mathfrak{Ig}_U^{(d_i)\sharp}$ over $\mathfrak{C}_U^{(b_i)\sharp}$ representable by a $(\zeta_1,\dotsc, \zeta_n)$-adic formal scheme which is finite \'etale over $\mathfrak C_U^{(b_i)\sharp}$. It is called a \emph{formal Igusa variety} of level $(d_i)$.

Similarly, $\mathfrak{Ig}_U^{(\infty_i)\sharp}$ parametrizing trivializations of the whole local $G_i$-shtukas is representable by a formal scheme, isomorphic to $\underset{\leftarrow}{\lim}~ \mathfrak{Ig}_U^{(d_i)\sharp}$.
\end{thm}
The foliation structure is then established by the following result.
\begin{thm}[\cite{Neupert}, Thm. 6.2.1]
There is a natural formally \'etale morphism of formal schemes over $\Spf A_c^{\sharp}$
$$\hat\pi_{(\infty_i)}:\prod_{i}\mathcal M_{b_{i,0}}^{\preceq \mu_i\sharp}\times_{\Spf A_c^{\sharp}} {\mathfrak{ Ig}_U^{(\infty_i)\sharp}}\rightarrow \mathfrak N_{U}^{[b_i]\sharp}.$$
\end{thm}
There are also variants of this result for the special fiber of the spaces, using bounds on the quasi-isogenies, and for the associated adic spaces in the sense of Huber.

The very rough idea of the construction of this morphism is the following: A point in the Igusa variety gives us a global shtuka, such that at each place $c_i$ the associated local $\G_i$-shtuka is the one defined by $b_{i,0}$. Then a point in the Rapoport-Zink space $\mathcal M_{b_{i,0}}^{\preceq \mu_i}$ gives us a modification by a quasi-isogeny of this local $\G_i$-shtuka, and we define a new global $\GG$-shtuka by keeping the old one away from the $c_i$ and replacing it in a formal neighborhood of $c_i$ by this new local $\G_i$-shtuka.

Applying the theorem to the cohomology of the moduli space of global $\GG$-shtukas, and its decomposition into Newton strata, one can express the cohomology of Newton strata, and therefore of the whole moduli space of global $\GG$-shtukas in terms of the cohomology of Rapoport-Zink moduli spaces and of Igusa varieties. For technical reasons the statement is only known over a finite extension $E'$ of $\mathbb{F}_q$ that can be described in terms of $\GG$ and $\mu$. 
 Denote by $K$ the function field of the curve $C$, by $\overline K$ an algebraic closure, and by $\Gamma_{E'}$ the absolute Galois group of $E'$.
\begin{thm}[\cite{Neupert}, Main Theorem 2]
Let $\nabla^{\mu}_n\mathcal{H}^1(C, \GG)$ be a moduli space of global $\GG$-shtukas, such that all connected components of $\nabla^{\mu}_n\mathcal{H}^1(C, \GG) \times_{C^n \setminus \Delta} \Spec E'\dbl \zeta_1,\dotsc, \zeta_n\dbr$ are proper over $\Spec E'\dbl \zeta_1,\dotsc, \zeta_n\dbr$. \\
Then there exists a canonical isomorphism between the virtual $\GG(\mathbb{A}^{c_i}_{\mathbb Q}) \times \Gamma_{E'}$-representations
$$\sum_j (-1)^j H^j_c\left(\nabla^{\mu}_n\mathcal{H}^1(C, \GG)\times \overline{K}, \mathbb{Q}_\ell\right)$$
and 
$$\sum_{[b_i]} \sum_{d, e, f} (-1)^{d+e+f} Tor_d^{\mathcal{H}(\prod J_{b_{i,0}})} \left(H^e_c\bigl(\prod \mathcal{M}_{b_{i,0}}^{\preceq \mu_i}\times \mathbb F, R\Psi_\eta^{\rm an}\mathbb{Q}_\ell ), \varinjlim_U \varinjlim_{d_i} H^f_c ({\rm Ig}^{(d_i)}_U\times \mathbb F, R\Psi_\eta^{\rm an}\mathbb{Q}_\ell \bigr)\right).$$
\end{thm}

\section{Period spaces and the period map}\label{secper}

In our context, period spaces are constructed as strictly $\BF_q\dpl\zeta\dpr$-analytic spaces in the sense of Berkovich. In the generality presented here, they are introduced and studied in \cite{hv3}. We allow more general bounds than those associated with minuscule coweights. Even if the bound is given by a coweight as in Example \ref{exbound}, the bound is in general a union of Schubert cells, and not a single one as in the minuscule case. For these two reasons the period spaces have to be defined as subspaces of an affine Grassmannian instead of a (classical) flag variety. To define them, we consider the affine Grassmannian $\Gr_{\G}^{B_{\rm dR}}$ of Remark \ref{remgrdr}. The space of Hodge-Pink $\G$-structures bounded by $\hat Z$ is defined as $\mathcal{H}_{\G,\hat Z}:=\hat Z_E$. It is a projective subscheme of $\Gr_G^{B_{\rm dR}}\widehat\times_{\BF_q\dpl\zeta\dpr} E$. Fix a local $\G$-shtuka $\underline{\mathbb G}_0$  over $\mathbb{F}$ and a trivialization $\underline{\mathbb G}_0\cong(L^+{\G}_\mathbb{F},b\sigma^*)$.

A $z$-isocrystal over $\mathbb{F}$ is a pair $(D,\tau_D)$ consisting of a finite-dimensional $\mathbb F\dpl z\dpr$-vector space $D$ and an isomorphism $\tau_D:\sigma^*D\rightarrow D$ of $\mathbb F\dpl z\dpr$-vector spaces. A Hodge-Pink structure on $(D,\tau_D)$ over a field extension $L$ of $\mathbb F\dpl\zeta\dpr$ is a free $L\dbl z-\zeta\dbr$-submodule $\mathfrak q_D\subset D\otimes_{\mathbb F\dpl z\dpr}L\dpl z-\zeta\dpr$ of full rank, compare \cite{HartlKim}, 5.  

\begin{definition}\label{defadm}
Let $L$ and $b$ be as above and let $\gamma\in\Gr_G^{B_{\rm dR}}(L)$. Let $\rho: G_{\mathbb F_q\dpl z\dpr}\rightarrow\GL_{n,\mathbb F_q\dpl z\dpr}$ be a representation and let $V$ be the representation space. Then we associate with $b$ (and $\rho$) the $z$-isocrystal
$(V\otimes_{\mathbb F_q\dpl z\dpr}\mathbb F\dpl z\dpr,\rho(\sigma^* b)\sigma^*)$  over $\mathbb F$.
With $\gamma$ we associate the Hodge-Pink structure over $L$ defined by
$$\mathfrak q_D(V)=\rho(\gamma)\cdot V\otimes_{\mathbb F_q\dpl z\dpr}L\dbl z-\zeta\dbr\subset D\otimes_{\mathbb F\dpl z\dpr}L\dpl z-\zeta\dpr.$$
Let $$\underline D_{b,\gamma}(V)=(V\otimes_{\mathbb F_q\dpl z\dpr}\mathbb F\dpl z\dpr,\rho(\sigma^* b)\sigma^*,\mathfrak q_D(V)).$$

Consider the Newton point $\nu_b\in X_*(T)$ as a homomorphism $\mathbb D_{\mathbb F\dpl z\dpr}\rightarrow G_{\mathbb F\dpl z\dpr}$. Recall the Kottwitz map $\kappa_G:G(\mathbb F\dpl z\dpr)\rightarrow \pi_1(G)_{\Gamma}$. We then define the \emph{Hodge degree} and the \emph{Newton degree} as
$$ t_H(\underline D_{b,\gamma}(V))=\rho_*(\kappa_G(\gamma))\qquad\text{ and }\qquad t_N(\underline D_{b,\gamma}(V))=\det_V\circ \rho\circ \nu_b.$$
Using $\pi_1(\GL(V))_\Gamma=\mathbb Z$ and $\Hom(\mathbb D_{\mathbb F\dpl z\dpr},\mathbb G_m)=\mathbb Q$ we can view these as rational numbers. (Compare \cite{hv3}, 4 for a definition in terms of $\underline D_{b,\gamma}(V)$ only.)

We call $\underline D_{b,\gamma}(V)$ \emph{weakly admissible} if the images of $[b]$ and $\gamma$ in $\pi_1(G)_{\Gamma,\mathbb Q}$ coincide and for all subobjects the Hodge degree is less or equal to the Newton degree. Finally, the pair $(b,\gamma)$ is weakly admissible if $\underline D_{b,\gamma}(V)$ is weakly admissible for a faithful representation $\rho$ of $G$. 

Let $\breve{\mathcal{H}}_{\G,\hat Z}:=\mathcal H_{\G,\hat Z}\widehat\times_{E}\breve E$ where $\mathcal H_{\G,\hat Z}=\hat Z_E$. The \emph{period space} $\breve{\mathcal{H}}_{\G,\hat Z,b}^{wa}$ is then defined as the set of all $\gamma$ in the associated $\breve E$-analytic space $\breve{\mathcal{H}}_{\G,\hat Z}^{\an}$ such that $(b,\gamma)$ is weakly admissible. 

To define the admissible locus $\breve{\mathcal{H}}_{\G,\hat Z,b}^{a}\subset\breve{\mathcal{H}}_{\G,\hat Z}^{\an}$, one associates with $\underline D$ a pair of $\sigma$-bundles $\underline{\mathcal E}(\underline D)=\underline{\mathcal E}_{b,\gamma}(V)$ and $\underline{\mathcal F}(\underline D)=\underline{\mathcal F}_{b,\gamma}(V)$. The \emph{admissible locus} is then defined as the subset over which $\underline{\mathcal F}_{b,\gamma}(V)$ has slope zero. Then $\breve{\mathcal{H}}_{\G,\hat Z,b}^{wa}$ and $\breve{\mathcal{H}}_{\G,\hat Z,b}^{a}$ are open paracompact strictly $\breve E$-analytic subspaces of $\breve{\mathcal{H}}_{\G,\hat Z}^{\an}$ (\cite{hv3}, Theorem 4.20).

The neutral admissible locus $\breve{\mathcal{H}}_{\G,\hat Z,b}^{na}\subset\breve{\mathcal{H}}_{\G,\hat Z,b}^{a}$ is defined by the condition that the images of $[b]$ and $\gamma$ in $\pi_1(G)_{\Gamma}$ coincide.
\end{definition}

We have the analog of the theorem of Colmez and Fontaine that ``admissible implies weakly admissible'', see \cite{HartlPSp},Theorem~2.5.3. In other words, we have an inclusion $$\breve{\mathcal{H}}_{\G,\hat{Z},b}^{a}\subseteq \breve{\mathcal{H}}_{\G,\hat{Z},b}^{wa}.$$ Furthermore, we have $\breve{\mathcal{H}}_{\G,\hat{Z},b}^{a}(L)=\breve{\mathcal{H}}_{\G,\hat{Z},b}^{wa}(L)$ for all finite field extensions $L/\breve E$. If $L$ is algebraically closed, weakly admissible does \emph{not} imply admissible in general.

For an $\mathbb F\dbl\zeta\dbr$-algebra $B$, complete and separated with respect to a bounded norm $|\cdot|: B\rightarrow [0,1]\subset \mathbb R$ with $0<|\zeta|<1$, we consider the $\mathbb F\dpl z\dpr$-algebra
$$
B\dbl z,z^{-1}\}=\{\sum_{i\in\mathbb Z} b_iz^i\mid b_i\in B,|b_i|\,|\zeta|^{ri}\rightarrow 0~(i\rightarrow-\infty) \text{ for all }r>0\}
$$
and the element $t_-=\prod_{i\in \mathbb N_0}\bigl(1-\tfrac{\zeta^{q^i}}{z}\bigr)\in B\dbl z,z^{-1}\}.$

We want to define the \emph{period morphism}
$$\breve\pi:\breve{\mathcal M}^\an\rightarrow\breve{\mathcal{H}}_{\G,\hat{Z},b}^a$$ as in \cite{hv3}, 6 where $\mathcal M$ is the Rapoport-Zink space associated with $\hat Z$ and $b$. Let $S$ be an affinoid, strictly $\breve E$-analytic space and let $S\rightarrow\breve{\mathcal M}^{\rm an}$ be a morphism of $\breve E$-analytic spaces. With it we have to associate a morphism $S\rightarrow\breve{\mathcal H}_{\G,\hat{Z}}^{\rm an}$. By construction of $\breve{\mathcal M}^\an$ a morphism $S\rightarrow \breve{\mathcal M}^\an$ is induced by a morphism from a quasi-compact admissible formal $\breve R_{\hat Z}$-scheme $\mathscr S$ with $\mathscr S^\an=S$ to $\breve{\mathcal M}$. The latter corresponds to some $(\underline{\mathcal G},\bar\delta)\in\breve{\mathcal M}(\mathscr S)$.
After an \'etale covering $\mathscr S'=\Spf B'\rightarrow\mathscr S$ of admissible formal $\breve R_{\hat Z}$-schemes there is a trivialization $\alpha:\underline{\mathcal G}_{\mathscr S'}\cong((L^+\G)_{\mathscr S'},A\sigma^*)$ for some $A\in \G(B'\dbl z\dbr[\tfrac{1}{z-\zeta}])$.

Denote by $\overline{\cdot}$ the reduction to $V(\zeta)$. Then the quasi-isogeny $\bar\delta$ induces via $\overline\alpha$ an element $\overline\Delta\in L\G(\overline{B}')$. One then shows that there is a uniquely determined element $ \Delta\in \G(B'\dbl z,z^{-1}\}[\tfrac{1}{t_{-}}])$ lifting $\overline \Delta$ with $\Delta A=b\,\sigma^*(\Delta)$ for $A$ as above and $b$ defining the Rapoport-Zink space. Consider the element
$$
\gamma  = \sigma^*(\Delta)A^{-1}\cdot \G(B'[\tfrac{1}{\zeta}]\dbl z-\zeta\dbr) \in \G(B'[\tfrac{1}{\zeta}]\dpl z-\zeta\dpr)/\G(B'[\tfrac{1}{\zeta}]\dbl z-\zeta\dbr).$$ The boundedness of $\underline{\mathcal G}$ then implies that $\gamma$ factors through $\breve{\mathcal H}_{\G,\hat{Z}}^{\rm an}$. It descends to a well-defined element of $\breve{\mathcal H}_{\G,\hat{Z}}^{\rm an}(S)$ which is the image under the period morphism.

\begin{prop}[\cite{hv3}, Prop.~6.9, 6.10] The period morphism factors through the open $\breve E$-analytic subspace $\breve{\mathcal H}_{G,\hat{Z},b}^a$ and induces an \'etale morphism
$$\breve\pi:\breve{\mathcal M}^\an\rightarrow\breve{\mathcal{H}}_{\G,\hat{Z},b}^a.$$
\end{prop}

To describe the image of $\breve \pi$ recall the tensor functors of \eqref{eqw0} and \eqref{eqforget}.
The construction of the $\sigma$-bundle $\underline{\mathcal{F}}_{b,\gamma}(V)$ of Definition \ref{defadm} has a generalization to a $\sigma$-bundle $\underline{\mathcal{F}}_b(V)$ over $\breve{\mathcal H}_{G,\hat{Z}}$ with fibers $\underline{\mathcal{F}}_{b,\gamma}(V)$. It further induces a canonical local system $\underline{\mathcal{V}}_b(V)$ of $\mathbb F_q\dpl z\dpr$-vector spaces on $\breve{\mathcal H}_{G,\hat{Z},b}^a$ with $\underline{\mathcal{V}}_b(V)_{\bar\gamma}=\bar\gamma^*\underline{\mathcal{F}}_{b}(V)^\tau$ where we take invariants under the isomorphism $\tau$ defining the $\sigma$-bundle. Let $\omega_{b,\bar\gamma}:\Rep_{\BF_q\dpl z\dpr}\Geta\rightarrow{\rm FMod}_{\BF_q\dpl z\dpr}$ be the fiber functor with $\omega_{b,\bar\gamma}(V):=\underline{\mathcal V}_b(V)_{\bar\gamma}=\bar\gamma^*\underline{\mathcal F}_b(V)^\tau$. Then, one of the main results of \cite{hv3} is the following description of the image of the period map. The condition that $G$ be unramified can also be replaced by another, more technical, but possible more general condition.

\begin{thm}\label{MainThm}
Assume that $G$ is unramified.
\begin{enumerate}
\item \label{MainThm_1}
The image $\breve\pi(\breve{\mathcal M}^{\rm an})$ of the period morphism is equal to the union of those connected components of $\breve{\mathcal H}_{\G,\hat{Z},b}^{a}$ on which there is an $\BF_q\dpl z\dpr$-rational isomorphism $\beta:\omega^{\circ}\rightarrow\omega_{b,\bar\gamma}$.
\item  \label{MainThm_2}
The rational dual Tate module $\check V_{\underline{\mathcal{G}}}$ of the universal local $\G$-shtuka $\underline{\mathcal{G}}$ over $\breve{\mathcal M}^{\rm an}$ descends to a tensor functor $\check V_{\underline{\mathcal{G}}}$ from $Rep_{\BF_q\dpl z\dpr}\Geta$ to the category of local systems of $\BF_q\dpl z\dpr$-vector spaces on $\breve\pi(\breve{\mathcal M}^{\rm an})$. It carries a canonical $J_b(\BF_q\dpl z\dpr)$-linearization and is canonically $J_b(\BF_q\dpl z\dpr)$-equivariantly isomorphic to the tensor functor $\underline{\mathcal V}_b$.
\item \label{MainThm_3}
The tower of strictly $\breve E$-analytic spaces $({\breve{\mathcal M}}^K)_{K\subset \G(\BF_q\dpl z\dpr)}$ is canonically isomorphic over $\breve\pi(\breve{\mathcal M}^{\rm an})$ in a Hecke and $J_b(\BF_q\dpl z\dpr)$-equivariant way to the tower of \'etale covering spaces of $\breve\pi(\breve{\mathcal M}^{\rm an})$ that is naturally associated with the tensor functor $\underline{\mathcal V}_b$.
\end{enumerate}
\end{thm}

\bigskip
  \footnotesize

  \textsc{Fakult\"at f\"ur Mathematik der
Technischen Universit\"at M\"unchen -- M11}

\textsc{Boltzmannstr. 3, 85748 Garching, Germany}

\texttt{viehmann@ma.tum.de}
\end{document}